\documentclass[a4paper,11pt,oneside]{amsart}

\usepackage{geometry}
\usepackage{graphicx, amsfonts, amssymb}
\usepackage{mathrsfs, amsmath, amsthm}
\usepackage{verbatim}
\usepackage{relsize}
\usepackage{enumerate}
\usepackage{hyperref}
\usepackage{version}
\usepackage{blindtext}
\usepackage{mathtools}

\newtheorem{theorem}{Theorem}[section]
\newtheorem*{thma}{Theorem A}

\newtheorem{Lemm}[theorem]{Lemma}
\newtheorem{prop}[theorem]{Proposition}

\newtheorem*{kahane}{Kahane's Inequality}
\theoremstyle{remark}
\newtheorem{remark}[theorem]{Remark}

\newcommand{\m}{\mathbb}

\DeclareMathOperator{\Hol}{Hol}
\DeclareMathOperator{\adj}{adj}
\hypersetup{
	colorlinks=true,
	linkcolor=blue,
	filecolor=magenta,      
	urlcolor=cyan,
	pdftitle={A Generalisation of Volterra Operator on Hardy spaces},
	pdfpagemode=FullScreen,
}

\title{On the boundedness of generalized integration operators on Hardy spaces}

\author{N. Chalmoukis}
\email{nikolaos.chalmoukis@unimib.it}
\address{Dipartimento di Matematica e Applicazioni, Universita degli studi di Milano Bicocca, via Roberto Cozzi, 55 20125, Milano, Italy}

\author{G. Nikolaidis}
\email{nikolaidg@math.auth.gr}
\address{Department of Mathematics, Aristotle University of Thessaloniki, 54124, Greece}

\thanks{This research project was supported by the Hellenic Foundation for Research and Innovation (H.F.R.I.) under the `2nd Call for H.F.R.I. Research Projects to support Faculty Members \& Researchers' (Project Number: 4662).}

\date{}

\begin{document}

\begin{abstract}

    We study the boundedness and compactness properties of the generalized integration operator $T_{g,a}$ when it acts between distinct Hardy spaces in the unit disc of the complex plane. This operator has been  introduced in \cite{chalmoukis2020generalized} by the first author in  connection to a theorem of Cohn about factorization of higher order derivatives of functions in Hardy spaces. We answer in the affirmative a conjecture stated in the same work, therefore giving a complete characterization of the class of symbols $g$ for which the operator is bounded from the Hardy space $H^p$ to $H^q, \, 0<p,q<\infty.$ 
\end{abstract}
\maketitle
\section{Introduction}
Let $\m{D}$ be the unit disc of the complex plane, $\m{T}$ be its boundary and  $\Hol(\m{D})$ the space of analytic functions defined in $\m{D}$. The classical Volterra operator is defined as
\begin{equation}\label{Voltera operator definition}
    Vf(z)=\int_{0}^{z}f(\zeta)\,d\zeta.
\end{equation}
For a fixed $g\in \Hol(\m{D})$,  we can define the following operator on $\Hol(\m{D})$ 
\begin{equation}\label{equation for Tg}
    T_gf(z)=\int_{0}^{z}f(\zeta)g'(\zeta)\,d\zeta\,,\qquad z\in\m{D}.
\end{equation}
The motivation to study integral operators such as $T_g$ comes partially from the fact that, as $g$ varies in $\Hol(\m{D})$, $T_g$ represents some significant classical operators. 
For instance, if $g(z)=z$, $T_g$ is the  Volterra operator $V$, while, when $g(z)=-\log(1-z)$, it coincides with the Ces\'aro summation operator. 
  Originally, the generalized Volterra operator was introduced and studied in the context of the Hardy spaces of the unit disc. The Hardy space $H^p$ is defined as the space of functions $f\in \Hol(\mathbb{D})$ such that

$$   \|f\|^p_{H^p}=\sup_{0\leq r<1}\int_{0}^{2\pi}|f(re^{it})|^p\,\frac{dt}{2\pi}<\infty\,.$$
	
	In particular, Ch. Pommerenke \cite{pommerenke1977schlichte}, characterized the symbols $g$ for which  $T_g$ is bounded on the Hilbert space $H^2$. 
The complete characterisation of symbols for which $T_g$ acts boundedly between different Hardy spaces was given in a series of papers (see \cite{aleman1995integral} and \cite{aleman2001integral}). 
Subsequently, the study of such operators on various spaces of analytic functions attracted a lot of attention (see \cite{aleman1997integration}, \cite{wu2006areabergman}, \cite{TanausuGalanopoulostentspaces}, \cite{pelaez2014weightedmemoir}). 
 
 	In this article we study a further generalization of the integral operator $T_g$.
An inspection of (\ref{equation for Tg}) shows that $T_g$ is a primitive of the first term of the derivative of the product $f g$. 
Applying the Leibniz rule of differentiation we get 
\begin{equation}\label{eq:Leibniz_dif}
 (T_gf)^{(n)}=\sum_{k=0}^{n-1}\binom{n-1}{k}f^{(k)} g^{(n-k)}.
\end{equation}

	Now, if we consider an arbitrary $n$-tuple $a=(a_0,a_1,\dots,a_{n-1})\in \mathbb{C}^{n}$, $a\neq \mathbf{0},$ we can define the following operator
 \begin{equation*}\label{equation defining Tga}
     T_{g,a}f=V^n\left(\sum_{k=0}^{n-1}a_kf^{(k)} g^{(n-k)}\right)\,,
 \end{equation*}
 where $V^n$ is the $n$-th iterate of the Volterra operator (\ref{Voltera operator definition}). 
 It is then clear by \eqref{eq:Leibniz_dif} that the generalized Volterra operator is 
 a particular instance of the opeartor $T_{g,a}.$ 

  The integral operator $T_{g,a}$ was introduced  by the first author in \cite{chalmoukis2020generalized}, in the context of Hardy spaces of the unit disc. 
Thereafter,
J. Du, S. Li, and D. Qu \cite{DuLiQu2022generalized} and X. Zhu \cite{xiaolingzhu2023generalizedtgnk}
studied the action of the operators
\begin{equation}\label{equation defining Tgnk}
T_g^{n,k}f= V^n(f^{(k)}g^{(n-k)})\,
\end{equation}
on weighted Bergman spaces and on $F(p,q,s)$ spaces 
respectively, which are specific examples of $T_{g,a}$ when $a$ is a standard unit vector of $\m{C}^n$. Recently, H. Arroussi et al. \cite{arroussi2024new} characterised the Sobolev-Carleson measures for  Bergman spaces and consequently characterised the space of symbols $g$ for which $T_{g,a}$ acts boundedly between Bergman spaces.

 The motivation for studying the operator $T_{g,a}$, beyond the fact that it generalizes the classical $T_g$, stems also from the connection of $T_{g,a}$ to a factorization theorem of holomorphic functions by W.Cohn \cite{cohn1999factorization} and a theorem of J. R\"atty\"a \cite{jounirattya2007lineareq} about higher order linear differential equations with holomorphic coefficients.
The interested reader is referred to \cite[Theorems 1.5 and 1.6]{chalmoukis2020generalized} for more details. 
 For completeness, we state the main result of \cite{chalmoukis2020generalized} regarding $T_{g,a}$. In order to do that, we recall the definitions of some spaces of analytic functions.

 Let $0< \alpha \leq 1.$ The analytic Lipschitz space $\Lambda_\alpha$ consists of $f\in \Hol(\m{D})$ which are continuous up to the boundary, and its boundary function $f(e^{i\theta})$ is H\"older continuous of order $\alpha$. An equivalent description, see \cite[Theorem 5.1]{duren1970theory}, is that $f\in\Lambda_{\alpha}$ iff
 
$$\|f\|_{\Lambda_{\alpha}}:=\sup_{z\in\m{D}}|f'(z)|(1-|z|^2)^{1-\alpha}<\infty.$$

The space $\lambda_{\alpha}$ consists of those functions in $f\in \Lambda_{\alpha}$, such that
$$\lim_{|z|\to 1}|f'(z)|(1-|z|^2)^{1-\alpha}=0.$$

A comprehensive introduction to this class can be found in \cite[Chapter 4]{duren1970theory}. Also, recall that the space of analytic functions of bounded mean oscillation, $BMOA$ consists of $f \in H^2$, such that 
\[ \sup_{a\in \mathbb{D}} \| f \circ \varphi_a - f(a) \|_{H^2} < + \infty,\]
where $\varphi_a(z) = \frac{z-a}{1-\overline{a}z}, z \in \mathbb{D}.$
Equivalent descriptions of $BMOA$ can be found for example in \cite{garnett2006bounded}. The following is the main theorem from \cite{chalmoukis2020generalized}. In what follows, $\m{Z}_{\geq0}=\m{N}\cup \{0\}.$

 \begin{thma}
 \hypertarget{Nikos result in Hardy}{L}et $0<p,q<\infty$, $n\in\m{N}$, $a=(a_0,\dots,a_{n-1})\in \m{C}^{n}$ with $a_0\neq 0$ and $g\in \Hol(\m{D})$.
     \begin{itemize}
         \item[(i)] When $0<p<q<\infty$, let $\ell= \max\{k\colon a_k\neq 0\}$ and $\alpha=\frac{1}{p}-\frac{1}{q}.$ Then the following is true.
         \begin{itemize}
             \item[(a)] If\, $\displaystyle{\kappa<\alpha\leq \kappa+1\leq n-\ell}$ for some $\kappa\in\m{Z}_{\geq0}$, then $T_{g,a}\colon H^p\rightarrow H^q$ is bounded if and only if
             \begin{equation*}\label{equation for p<q in Tga in Hardy}
                g^{(\kappa)}\in\Lambda_{\alpha-\kappa}\,.
             \end{equation*}
             \item[(b)] If $\displaystyle{\alpha> n-\ell}$ and $T_{g,a}\colon H^p\rightarrow H^q$ is bounded, then $T_{g,a}$ is the zero operator.
         \end{itemize}
         \item[(ii)] When $p=q$, then $T_{g,a}\colon H^p\rightarrow H^p$  is bounded if and only if
            $$g\in {BMOA}\,.$$
         \item[(iii)] When $0<q<p<\infty$ and $g\in H^{\frac{pq}{p-q}}$, then $T_{g,a}\colon H^p\rightarrow H^q$ is bounded. When $n=2$, $a=(1,0)$ and assuming that $T_{g,a}\colon H^p\to H^q$ is bounded, then
         $$g\in H^{\frac{pq}{p-q}}\,.$$
     \end{itemize}
 \end{thma}
 As it is clear, the third part of Theorem \hyperlink{Nikos result in Hardy}{A} does not give a complete characterisation for the case $0<q<p<\infty$. Moreover, this result does not offer us any information when $a_0=0$. The main scope of this article is to answer completely these questions.
\begin{theorem}\label{characterisation for 0<q<p<+infty}
     Let $0<q<p<\infty$, $n\in\m{N}$, $g\in \Hol(\m{D})$ and $a=(a_0,\dots,a_{n-1})\in \mathbb{C}^{n}$, $a\neq \mathbf{0}$.
     \begin{itemize}
         \item[(i)] If $a_0\neq 0$, then $T_{g,a}\colon H^p\to H^q$ is bounded, if and only if $g\in H^{\frac{pq}{p-q}}.$
         \item[(ii)] If $a_0=0$, then $T_{g,a}\colon H^p\to H^q$ is bounded, if and only if $g\in BT^{\frac{pq}{p-q}}.$
     \end{itemize}
\end{theorem}

 The space $BT^p$, $0<p<\infty$, sometimes called the Bloch tent space, consists of $f\in\Hol(\m{D})$ such that 
\begin{equation}\label{Blochtentseminorm}
\|f\|^p_{BT^p}=\int_{\m{T}}\left(\sup_{z\in\Gamma_M(\zeta)}|f'(z)|(1-|z|^2)\right)^p|d\zeta|<\infty.
\end{equation}
where, $|d\zeta|$ is the normalized arc-length measure on $\mathbb{T}$ and $\Gamma_M(\zeta)$ is the Stolz angle of aperture $M$ centred at $\zeta$, i.e.,
$$\Gamma_M(\zeta)=\biggl\{z\in\m{D}\colon |1-z\overline{\zeta}|< \frac{M}{2}(1-|z|^2)\biggr\}\,,\qquad M>1, \,\zeta\in\m{T}.$$
 This space of analytic functions has been studied only recently and there is limited literature describing its properties. 
In fact, the first to study the basic properties of $BT^p$ was A. Per\"al\"a in \cite{PeralaBlochtent}. In addition to other results, he proved that this space can be identified using also higher order derivatives of functions when $p>1$ \cite[Theorem 4]{PeralaBlochtent}. 
Additionally, properties such as growth properties of functions in $BT^p$  have been recently studied in \cite{chen2023closuresblochtent}. 
For our purposes, we prove that we can express the seminorm $\|\cdot \|_{BT^p}$ using higher order derivatives of functions, extending the result of \cite{PeralaBlochtent} for $0<p\leq 1$.

 Finally, we finish the study of the integral operator $T_{g,a}$ by characterising the space of symbols $g$ for which the operator $T_{g,a}$ acts compactly between Hardy spaces when $a_0\neq 0$. Incorporating the  known result \cite[Theorem 1.1]{chalmoukis2020generalized} we arrive at the following theorem.

\begin{theorem}\label{compactness characterisation on Hardy spaces}
Let $0<p,q<\infty$\,, $n\in\m{N}$\,, $a=(a_0,\dots,a_{n-1})\in \m{C}^{n}$ with $a_0\neq 0$ and $g\in \Hol(\m{D})$.
     \begin{itemize}
         \item[(i)] When $0<p<q<\infty$, let $\ell= \max\{k\colon a_k\neq 0\}$  and $\alpha=\frac{1}{p}-\frac{1}{q}.$ Then the following conditions hold.
         \begin{itemize}
             \item[(a)] If\, $\kappa<\alpha\leq \kappa+1< n-\ell$ for some $\kappa\in\m{Z}_{\geq0}$, then $T_{g,a}\colon H^p\rightarrow H^q$ is compact if and only if
             \begin{equation}\label{equation for p<q in Tga compact Hardy}
                 g^{(\kappa)}\in\lambda_{\alpha-\kappa}\,.
             \end{equation}
             \item[(b)] If $\alpha= n-\ell$ and $T_{g,a}\colon H^p\rightarrow H^q$ is compact, then $T_{g,a}$ is the zero operator.
         \end{itemize}
         \item[(ii)] When $0<q<p<\infty$, $T_{g,a}\colon H^p\rightarrow H^q \text{ is compact}$ whenever it is bounded, i.e 
         $$g\in H^{\frac{pq}{p-q}}.$$
     \end{itemize}
\end{theorem}

Our techniques fall short of characterizing the compactness of $T_{g,a}: H^p \to H^q$ when $q<p$ and $a_0=0.$

  As it is customary, for real valued functions $A,\, B$, we write $A\lesssim B$, if there exists a positive constant $C$ (which may be different in each occurrence) independent of the arguments of $A,\,B$ such that $ A\leq CB$. The notation $A\gtrsim B$ can be understood in an analogous manner. If both $A\lesssim B$ and $A\gtrsim B$ hold simultaneously, then we write $A\cong B$.

  \section{Preliminaries}
  \subsection{Hyperbolic geometry}We recall some elementary facts from the geometry of the Poincar\'e disc. Let
$$dA(z)=\frac{dxdy}{\pi},\qquad z=x+iy$$
be the normalised Lebesgue area measure on $\m{D}$. We use the notation
 $$\rho(z,w)=|\varphi_z(w)|=\left|\frac{z-w}{1-\overline{z}w}\right|,\qquad z,w\in\m{D}$$
 for the {\it pseudohyperbolic} metric of $\m{D}$. Moreover,
 $$\beta(z,w)=\frac{1}{2}\log\frac{1+\rho(z,w)}{1-\rho(z,w)},\qquad z,w\in\m{D}$$
 is the {\it hyperbolic} metric of $\m{D}$.
  The set $\displaystyle{D(z,r)= \{w\in\m{D}\colon \beta(z,w)<r\}}$ is the hyperbolic disc, centered at $z$ with radius $r>0$. For $f\in \Hol(\m{D})$, by applying the Cauchy integral formula and subsequently a sub-mean inequality for hyperbolic discs \cite[Lemma 13, p. 66]{duren2004bergman}, we have that, for every $n\in\m{N},$ $0<p<\infty$ and $r>0$ , there exists $C=C(n,r)>0$ such that for all  $z\in\m{D}$ and $f\in \Hol(\mathbb{D})$, 
  \begin{equation}\label{Luecking hyperbolic estimate}
     |f^{(n)}(z)|^p\leq \frac{C}{(1-|z|^2)^{2+np}}\int_{D(z,r)}|f(w)|^p\,dA(w).
 \end{equation}
\indent 

It is well known, see \cite[Proposition 4.5]{zhu2007operator}, that given $z\in\m{D}$ and $r>0$,
\begin{align}
    |1-\overline{w}a|&\cong |1-\overline{w}z|\,,\qquad w\in \m{D},a\in D(z,r),\nonumber\\
    1-|w|^2&\cong 1-|z|^2\,,\qquad w\in D(z,r),\label{geometric estimates}\\
    A(D(z,r))& \cong (1-|z|^2)^2.\nonumber
\end{align}

Additionally, the following result, which connects the hyperbolic distance and the region $\Gamma_M(\zeta)$, is crucial in our work.
 \begin{Lemm}\label{lemma of covering hyperbolic disc the stoltz angle}
    Let $M>0$ and $r>0$. If $M^*=(M+1)e^{2r}-1>M$ then
$$\bigcup_{z\in\Gamma_M(\zeta)}D(z,r)\subset \Gamma_{M^*}(\zeta).$$
\end{Lemm}

    For a proof of this result, see \cite[Lemma 2.3, p. 992.]{wu2006areabergman}.

A sequence $\{z_\lambda\}_\lambda\subset \m{D}$ is called $r$ - {\it hyperbolically  separated} if there exists  a constant $r>0$ such that $\beta(z_k,z_\lambda)\geq r $ for $k\neq\lambda$, while is said to be an $(R, r)$-lattice in the hyperbolic distance, for $R > r > 0$, if it is $2r$ - separated and  
\[ \mathbb{D} = \bigcup_k D(z_k, R). \]

 When the particular constants are not important we will say simply lattice. Hyperbolically separated sequences have the following useful ``finite overlapping'' property.

\begin{Lemm}\label{finite covering Lemma}
    Let $r>0$ and $Z=\{z_\lambda\}_\lambda$ be a hyperbolically separated sequence.
    There exists a constant $N=N(r)>0$ such that for every $z\in\m{D}$ there exist at most $N$ hyperbolic discs $D(z_\lambda,r)$ such that $z\in D(z_\lambda,r).$  
\end{Lemm}
A proof can be found in \cite[Lemma 4.8]{zhu2007operator}.

  \subsection{Function spaces} 

  In the following section, we recall some known facts about the spaces of functions we consider, while also prove some auxiliary results we will need later.

Given $0<p<\infty$, the point evaluation functionals of the derivatives are bounded in the Hardy spaces $H^p$ \cite[Lemma p. 36]{duren1970theory}. 

That is, for every $z\in\m{D}$, there exists a constant $C=C(n,p)>0$, such that
\begin{equation}\label{equation for pointevaluations}
    |f^{(n)}(z)|\leq \frac{C}{(1-|z|^2)^{n+\frac{1}{p}}}\|f\|_{H^p}\,,\qquad\forall f\in H^p,\,\, n\in\m{N}.
\end{equation}
 Since $T_{g,a}$ is an integral operator, we desire a way to connect the $H^p$ norm with a quantity which involves the derivatives. By a well known result of C. Fefferman and E. Stein \cite{fefferman1972h} and its extension to higher order derivatives by P. Ahern and J. Bruna \cite[Theorem 4.2]{ahern1988maximal}, we have that
\begin{equation}\label{tent norm of Hardy derivatives}
     \|f\|_{H^p}^p\cong \sum_{k=0}^{n-1}|f^{(k)}(0)|^p+\int_{\m{T}}\left(\int_{\Gamma_M(\zeta)}|f^{(n)}(z)|^2(1-|z|^2)^{2n-2}\,dA(z)\right)^{p/2}|d\zeta|.
\end{equation}

 Let $0<p,q<\infty$. The {\it tent spaces} $T_q^p$ consist of measurable functions $\varphi$ defined on $\m{D}$ such that
$$\|\varphi\|_{T_q^p}= \left(\int_{\m{T}}\left(\int_{\Gamma_M(\zeta)}|\varphi(z)|^qdA(z)\right)^{\frac{p}{q}}|d\zeta|\right)^{1/p}<\infty\,.$$
When $q=\infty$, we define the tent space in terms of the non-tangential maximal function. Namely, the space $T_\infty^p$ consist of measurable functions $\varphi$ defined on $\m{D}$ such that
$$\|\varphi\|_{T_\infty^p}=\left(\int_{\m{T}}N_M(\varphi)(\zeta)^p\,|d\zeta|\right)^{1/p}<\infty,$$
where 
$$N_M(\varphi)(\zeta)=\sup_{z\in \Gamma_M(\zeta)}|\varphi(z)|$$
is the non-tangential maximal function.
   \begin{remark}\label{remark for well defined bloch tent} While the quantities $N_M(\varphi)$ and $\int_{\Gamma_M(\zeta)}|\varphi(z)|^q\,dA(z)$ both depend on the value of $M$, the spaces $T_p^q$ and $T_\infty^p$ do not. This follows for $T_{\infty}^p$ by \cite[Lemma 1, p. 166]{fefferman1972h}, while for the other values of $p,q$, it follows from \cite[Proposition 1]{LueckingderivativesHardy}. Consequently, from now on, we drop the subscript $M$, where it does not play any role in the arguments.
   \end{remark}
  Proceeding further, we focus our attention to the Bloch tent space.  Even though \eqref{Blochtentseminorm} only induces a semi norm on $BT^p$, it suffices for our needs, since $T_{g+c,a}=T_{g,a}$ for $c$ a constant. 
    
It follows immediately that the Bloch space $ \mathcal{B} = BT^\infty$ is contained in $BT^p, p>0$.
In this article, we need to provide a way to identify $f\in BT^p$ functions by a quantity involving higher order derivatives of $f$. We mention that in \cite[Theorem 4]{PeralaBlochtent}, Per\"al\"a proves the same result for $p>1$. However, as his method include duality arguments, it does not extend to $0<p\leq 1.$ 
\begin{prop}\label{Blochtentderivatives}
    Let $n\in\m{N}$, $0<p<\infty$ and $f\in \Hol(\m{D})$. Then $H^p \subseteq  BT^p$ and $ f \in BT^p$ if and only if 

$$\int_{\m{T}}\left(\sup_{z\in\Gamma(\zeta)}|f^{(n)}(z)|(1-|z|^2)^n\right)^p|d\zeta|<\infty.$$
Moreover,
$$\|f\|_{BT^p}^p\cong \sum_{k=1}^{n-1}|f^{(k)}(0)|+\int_{\m{T}}\left(\sup_{z\in\Gamma(\zeta)}|f^{(n)}(z)|(1-|z|^2)^n\right)^p|d\zeta|\,.$$
\end{prop}
\begin{proof}
    Let $M>1$ and fix a $\zeta\in\m{T}$. By the means of \eqref{Luecking hyperbolic estimate}, we have that, for $z\in\Gamma_M(\zeta)$,
    $$|f^{(n)}(z)|\lesssim \frac{1}{(1-|z|^2)^{n+1}}\int_{D(z,r)}|f'(w)|dA(w).$$
    Hence, using \eqref{geometric estimates}, we conclude that,
    \begin{align*}
        (1-|z|^2)^{n}|f^{(n)}(z)|&\lesssim \frac{1}{(1-|z|^2)}\int_{D(z,r)}|f'(w)|\,dA(w)\\
        & \lesssim \frac{1}{(1-|z|^2)^2}\int_{D(z,r)}(1-|w|^2)|f'(w)|dA(w)\\
        & \leq  \frac{A(D(z,r))}{(1-|z|^2)^2}\sup_{w\in D(z,r)}|f'(w)|(1-|w|^2) \\
        & \cong \sup_{w\in D(z,r)}|f'(w)|(1-|w|^2).
    \end{align*}
    For that fixed $r>0$, Lemma \ref{lemma of covering hyperbolic disc the stoltz angle} implies that there exists a $M^*>M>1$, such that $D(z,r)\subset \Gamma_{M^*}(\zeta)$ for all $z\in\Gamma_M(\zeta).$ Hence, for every $z\in\Gamma_M(\zeta)$, we conclude that
    \begin{align*}
        (1-|z|^2)^{n}|f^{(n)}(z)|&\lesssim  \sup_{w\in D(z,r)}|f'(w)|(1-|w|^2)\\
        & \leq \sup_{w\in \Gamma_{M^*}(\zeta)} |f'(w)|(1-|w|^2).
    \end{align*}

    Hence, the independence of the definition of the tent Bloch space from the aperture of the Stolz angle, allows us to conclude that

    $$\int_{\m{T}}\left(\sup_{z\in\Gamma_{M}(\zeta)}|f^{(n)}(z)|(1-|z|^2)^n\right)^p|d\zeta|\lesssim \|f\|_{BT^p}^p.$$
    Let now $f\in H^p$. By \cite[Theorem 3.1]{garnett2006bounded}, we have that $f\in T_{\infty}^p$. 

    Similarly, we have that
    \begin{align*}
        (1-|z|^2)|f'(z)|&\lesssim \sup_{w\in D(z,r)}|f(w)|.
    \end{align*}
    
    Hence, as above, there exists $M^*>M$ such that for every $z\in\Gamma_M(\zeta)$, it holds
    \begin{align*}
        (1-|z|^2)|f'(z)|&\lesssim  \sup_{w\in \Gamma_{M^*}(\zeta)} |f(w)|.
    \end{align*}
    Therefore, we conclude that for every $f\in H^p$, we have that $f\in BT^p$, proving the first part of this proposition. 
     For the final conclusion, fix again an $M>1.$ Without the loss of generality, we assume that 
    $$f(0)=f'(0)=\dots=f^{(n-1)}(0)=0.$$
    
    Observe that for every $\zeta\in\m{T}$, the Stoltz region $\Gamma_M(\zeta)$ is a convex set, containing the point $0$. Consequently, if $z\in\Gamma_M(\zeta),$ the line segment which connects $0$ with $z$ denoted by $[0,z]$, lies inside $\Gamma_M(\zeta).$ If $z=|z|e^{i\theta}$, then
    \begin{align*}
        |f^{(n-1)}(z)|& \leq \int_{0}^{|z|}|f^{(n)}(te^{i\theta})|\,dt\\
        & \leq \sup_{\xi\in[0,z]}|f^{(n)}(\xi)|(1-|\xi|)^n \int_{0}^{|z|}\frac{1}{(1-t)^n}\,dt\\
        & \leq  \sup_{\xi\in\Gamma_M(\zeta)}|f^{(n)}(\xi)|(1-|\xi|^2)^n  \cdot \frac{1}{(1-|z|)^{n-1}}.
    \end{align*}
    Thus,
    $$\int_{\m{T}}\left(\sup_{z\in\Gamma_M(\zeta)}|f^{(n-1)}(z)|(1-|z|^2)^{n-1}\right)^p|d\zeta|\leq \int_{\m{T}}\left(\sup_{z\in\Gamma_M(\zeta)}|f^{(n)}(z)|(1-|z|^2)^{n}\right)^p|d\zeta|.$$
    An induction on $n$ gives then desired result.
\end{proof}

Finally, we turn our attention to analytic Lipschitz classes. Once more, we shall need the equivalent description of the functions in these spaces using higher order derivatives. We refer the reader to \cite{zhu1993bloch} for a proof.
\begin{prop}\label{Proposition Lipschitz class derivatives}
Let $0<\alpha\leq 1$, $n\in\m{N}$. The following  hold;
\begin{itemize}
\item[(i)] A holomorphic function $f$ belongs to $\Lambda_\alpha$ if and only 
 \begin{equation*}
     \sup_{z\in\m{D}}(1-|z|^2)^{n-\alpha}|f^{(n)}(z)|<\infty\,.
 \end{equation*}
\item[(ii)] A holomorphic function $f$ belongs to $ \lambda_{\alpha}$ if and only 
 \begin{equation*}
     \lim_{|z|\to1}(1-|z|^2)^{n-\alpha}|f^{(n)}(z)|=0\,.
 \end{equation*}
 \end{itemize}
  \end{prop}

\subsection{Tent sequence spaces}
     
 The discrete analogue of tent spaces are the {\it tent sequence spaces}. Let $Z=\{z_\lambda\}_\lambda$, $\lambda\in \m{Z}_{\geq 0},$ be a lattice and $0<p,q<\infty$. We define the {\it tent sequence space} $T_q^p(Z)$, consisting of complex sequences $\{c_\lambda\}_\lambda$ such that
    $$\|\{c_\lambda\}\|_{T_q^p(Z)}=\left(\int_{\m{T}}\left(\sum_{\lambda\colon z_\lambda\in \Gamma(\zeta)}|c_\lambda|^q\right)^{p/q}|d\zeta|\right)^{1/p}<\infty.$$
    If $q=\infty$, we define the space $T_{\infty}^p(Z)$ as the sequence space of $\{c_\lambda\}_\lambda$ satisfying
    $$\|\{c_\lambda\}\|_{T_{\infty}^p(Z)}=\left(\int_{\m{T}}\left(\sup_{\lambda\colon z_\lambda\in \Gamma(\zeta)}|c_\lambda|\right)^{p}|d\zeta|\right)^{1/p}<\infty.$$ 
    
Such spaces have been used to study derivative embedding problems in Hardy spaces \cite{LueckingderivativesHardy}. For our purposes, we mention the factorization and the duality properties of tent sequence spaces that we need.
\begin{Lemm}\label{factorization of tent sequences}
    Let $Z=\{z_\lambda\}_\lambda$ be a lattice and $0<p\leq p_1,p_2\leq \infty$ and $0<q\leq q_1,q_2\leq \infty$ satisfying  
    $$\frac{1}{p}=\frac{1}{p_1}+\frac{1}{p_2}, \qquad \frac{1}{q}=\frac{1}{q_1}+\frac{1}{q_2}.$$
    If $\{c_{\lambda}\}_{\lambda}\in T_{p_1}^{q_1}(Z)$ and $\{t_\lambda\}_{\lambda}\in T_{p_2}^{q_2}(Z)$, then $\{c_\lambda t_\lambda\}_{\lambda}\in T_{p}^{q}(Z)$ with
    $$\|\{c_\lambda t_\lambda\}\|_{T_{p}^{q}(Z)}\lesssim \|\{c_\lambda\}\|_{T_{p_1}^{q_1}(Z)}\cdot \|\{t_\lambda\}\|_{T_{p_2}^{q_2}(Z)}.$$
    Conversely, if $\{k_\lambda\}_{\lambda}\in T_p^q(Z)$, then there exist sequences $\{c_{\lambda}\}_{\lambda}\in T_{p_1}^{q_1}(Z)$ and $\{t_\lambda\}_{\lambda}\in T_{p_2}^{q_2}(Z)$ such that $k_\lambda=c_\lambda t_\lambda$ and
    $$\|\{c_\lambda\}\|_{T_{p_1}^{q_1}(Z)}\cdot \|\{t_\lambda\}\|_{T_{p_2}^{q_2}(Z)}\lesssim \|\{k_\lambda \}\|_{T_{p}^{q}(Z)}\,.$$
    
\end{Lemm}

\begin{Lemm}\label{duality of tent sequence spaces}
    Let $1<q<\infty$. The dual of $T_1^q(Z)$ can be identified with the space $T_{\infty}^{q'}(Z)$, where $q'=\frac{q}{q-1}$, under the pairing
    $$\langle c_\lambda,\mu_\lambda\rangle = \sum_{\lambda}c_\lambda \mu_\lambda (1-|z_\lambda|),$$
    where $\{c_\lambda\}_\lambda\in T_1^q(Z)$ and $\{\mu_\lambda\}_\lambda\in T_{\infty}^{q'}(Z).$ In particular, 
    $$\|\{\mu_{\lambda}\}\|_{T_{\infty}^{q'}(Z)}\cong \sup\biggl\{\Bigl|\sum_{\lambda}c_\lambda \mu_{\lambda}(1-|z_\lambda|)\Bigr|\colon \quad \| \{c_\lambda\} \|_{T_1^q(Z)}=1\biggr\}\,.$$
\end{Lemm}
See \cite[Proposition 6]{miihkinen2020volterra} for a proof of Lemma \ref{factorization of tent sequences} and \cite[Proposition 2]{LueckingderivativesHardy}, (see also \cite[Lemma 6]{arsenovic1999embedding}) for a proof of Lemma \ref{duality of tent sequence spaces}.
 The bridge between the boundedness of $T_{g,a}$ and tent sequence space is the following result, which is a slight modification of \cite[Lemma 3]{LueckingderivativesHardy}, so that we omit its proof. We mention also that this result is proved in \cite[Proposition A]{lv2021pautent} for tent spaces in the unit ball of $\m{C}^n$.
\begin{Lemm}\label{Operator Lemma}
    Let $0<p<\infty$, $b>\max\{1,2/p\}$, $j\in \m{Z}_{\geq0}$ and $Z=\{z_\lambda\}_\lambda$ be a lattice. Then the function
    $$S_j[\{c_\lambda\}](z)=\sum_{\lambda}c_\lambda \left(\frac{1-|z_\lambda|^2}{1-|z_\lambda|^j\overline{z_\lambda}z}\right)^b$$
    belongs to $H^p$, whenever $\{c_\lambda\}_\lambda\in T_2^p(Z)$ and 
    $$\|S_j[\{c_\lambda\}]\|_{H^p}\lesssim\|\{c_\lambda\}\|_{T_2^p(Z)}\,. $$ 
    
\end{Lemm}
\subsection{Some further lemmas}

Finally, we present the auxiliary lemmas we shall use in later sections.
The first inequality is the well known Banach valued variant of Khinchine's inequality. 
\begin{kahane}
    \hypertarget{kahane}{D}efine the Rademacher functions $r_\lambda$ by
    \begin{align*}
        r_0(t)&=\begin{cases}
        1,& 0\leq t-[t]<\frac{1}{2}\\
        -1,& \frac{1}{2}\leq t-[t]<1.
    \end{cases}\\
    r_\lambda(t)&= r_0(2^{\lambda}t)\qquad \lambda\geq 1.
    \end{align*}
Also, let $(X,\|\,\|)$ be a Banach space and $0<p<q < \infty$. There exists positive constant $a=a(p,q), b=b(p,q)$ such that for all $m\in\m{N}$ and $x_1,x_2,\dots,x_m\in X$
$$a\left(\int_{0}^{1}\left\|\sum_{\lambda=1}^{m}r_\lambda(t)c_\lambda\right\|^pdt\right)^{1/p}\leq \left(\int_{0}^{1}\left\|\sum_{\lambda=1}^{m}r_\lambda(t)c_\lambda\right\|^qdt\right)^{1/q}\leq b\left(\int_{0}^{1}\left\|\sum_{\lambda=1}^{m}r_\lambda(t)c_\lambda\right\|^pdt\right)^{1/p} .$$
 \end{kahane}
 The next lemma is a simple linear algebra calculation that we require for the necessity part in Theorem \ref{characterisation for 0<q<p<+infty}.

 \begin{Lemm}\label{algebraic lhmma}
    Let $0<p<\infty$, $a=(a_0,\dots,a_{n-1})\in\m{C}^n$ and assume that $f_0,f_1,\dots,f_{n-1}$ are complex valued functions on the unit disc. Given the system of linear equations
    $$D_j(z) = \sum_{k=0}^{n-1}|z|^{jk}f_k(z)\frac{(1-|z|^2)^n}{(1-|z|^{j+2})^k}\qquad j=0,\dots,n-1\,,$$
    then for each $0\leq k\leq n-1$,
    $$f_k(z)(1-|z|)^{n-k}=\sum_{j=0}^{n-1}b_{jk}(z) D_j(z),\qquad 0<
    |z|<1,$$
    where $b_{jk}$ are bounded when $\frac{1}{2}<|z|<1.$
\end{Lemm}
\begin{proof}
    The proof of the lemma amounts to solving a system of linear equations. 
  
   The original assumption is equivalently translated in the following matrix equation;
  $$  \begin{bmatrix}
        1& 1&\dots &1\\
        1& \frac{|z|(1+|z|)}{1+|z|+|z|^2}&\dots &\left(\frac{|z|(1+|z|)}{1+|z|+|z|^2}\right)^{n-1}\\
        \vdots&\vdots&\ddots&\vdots\\
        1& \frac{|z|^{n-1}(1+|z|)}{\sum_{k=0}^{n}|z|^k}&\dots &\left(\frac{|z|^{n-1}(1+|z|)}{\sum_{k=0}^{n}|z|^k}\right)^{n-1}
    \end{bmatrix} 
    \begin{bmatrix*}[l]
        &f_0(z)(1-|z|)^n\\
        &f_1(z)(1-|z|)^{n-1}\\
        &\vdots\\
        &f_{n-1}(z)(1-|z|)
    \end{bmatrix*}=
    \begin{bmatrix*}[c]
        D_0(z)\\
        D_1(z)\\
        \vdots\\
        D_{n-1}(z)
    \end{bmatrix*}.$$
    The matrix on the left hand side is a Vandermonde matrix, $V=\{x_i^j\}_{i,j=0}^{n-1}$, with
    $$x_i = x_i(z) =\frac{|z|^i(1+|z|)}{\sum_{k=0}^{i+1}|z|^k}.$$
    The determinant of $V$ is given by
$$\det(V)=\prod_{0\leq i<\lambda\leq n-1}(x_i-x_\lambda)\neq 0, $$
when $0<|z|<1.$ The inverse matrix of $V$, is given by the equation $V^{-1} = \det(V)^{-1} \adj(V)$, where $\adj(V)$ is the adjoint matrix. 
We m.ust show that the elements of $V^{-1}$, as functions of $z$, stay bounded for $1/2\leq |z|<1.$ As $\operatorname{adj}(V)$ contains elements which are polynomials of the elements of $V$, then the elements of $\operatorname{adj}(V)$ are bounded in $\frac{1}{2}\leq |z|<1$, as $x_i(z)$ are bounded. To finish the proof, it suffices to show a lower bound for $\det(V)$. In particular, we compute that for  $j<\nu$,
\begin{align*}
    |x_\nu(z)-x_j(z)|& =\left|\frac{|z|^{\nu}(1+|z|)}{\sum_{k=0}^{\nu+1}|z|^k}-\frac{|z|^j(1+|z|)}{\sum_{k=0}^{j+1}|z|^k}\right|\\
    & =|z|^j(1+|z|)\left|\frac{|z|^{\nu-j}\sum_{k=0}^{j+1}|z|^k-\sum_{k=0}^{\nu+1}|z|^k}{\sum_{k=0}^{j+1}|z|^k\sum_{k=0}^{\nu+1}|z|^k}\right|\\
    & \geq \left(\frac{1}{2}\right)^{n-1}\frac{\sum_{k=0}^{\nu+1}|z|^k-|z|^{\nu-j}\sum_{k=0}^{j+1}|z|^k}{(j+2)(\nu+2)}\\
    & \geq \left(\frac{1}{2}\right)^{n-1}\frac{1}{(j+2)(n-1)}\sum_{k=j+2}^{\nu+1}|z|^k\\
    & \geq \left(\frac{1}{2}\right)^{n+j+1}\frac{1}{(j+2)(n-1)}\,.
\end{align*}
\end{proof}
Finally, let us recall \cite[Lemma 2.3]{chalmoukis2020generalized} which we state for completeness of the presentation. Here and subsequently,  
$$(\gamma)_0=1\quad\text{and}\quad (\gamma)_k=(\gamma)_{k-1}(\gamma+k-1)\qquad k\geq 1.$$
\begin{Lemm}\label{Nikos Lemma for test functions}
    Suppose that $f_0, f_1,\dots, f_{n-1}$ are complex valued functions on the unit disc and $\gamma$ be sufficiently large. If for any $\{z_m\}_m\subset \m{D}$ such that $|z_m|\to 1$ as $m\to\infty$, we have that
$$\lim_{m\to \infty}\left|\sum_{k=0}^{n-1}f_k(z_m) (\gamma)_k\right|=0\,,$$
then,
$$\lim_{m\to\infty}|f_k(z_m)|=0\qquad \forall\,\, 0\leq k\leq n-1\,.$$ 
\end{Lemm}

 \section{Proof of main results}
  To prove Theorem \ref{characterisation for 0<q<p<+infty}, we study first the operators $T_g^{n,k}$, we defined in \eqref{equation defining Tgnk}. Clearly, the operators $T_{g,a}$ are linear combinations of operators $T^{n,k}_g$.   The following result extends the previous result of \cite[Theorem 1.3]{chalmoukis2020generalized}. For the proofs that will follow in this section, set $$dA_n(z)=(1-|z|^2)^{2n-2}dA(z)\qquad n\in\m{N}\,.$$

 \begin{prop}\label{Extension of Nikos tgnk result}
     Let $0<q<p<\infty$, $n\in\m{N}$ and $1\leq k\leq n-1$. If $g\in BT^{\frac{pq}{p-q}}$, then $T_g^{n,k}\colon H^p\to H^q$ is bounded and $\|T_g^{n,k}\|_{H^p\to H^q}\lesssim \|g\|_{BT^{\frac{pq}{p-q}}}$.
 \end{prop}
\begin{proof}
     Let $f\in H^p$. By estimating  the Hardy norms of $T_{g}^{n,k}f$ and $f$ using \eqref{tent norm of Hardy derivatives}, applying H\"older's inequality, and  estimating the $BT^p$  seminorm of $g$ by the means of Proposition \ref{Blochtentderivatives}, we conclude that
    \begin{align*}
	    \|T_{g}^{n,k}f\|_{H^q}^q&\cong \int_{\m{T}}\left(\int_{\Gamma(\zeta)}|(T_{g}^{n,k}f)^{(n)}(z)|^2dA_n(z)\right)^{q/2}|d\zeta|\\
        &=\int_{\m{T}}\left(\int_{\Gamma(\zeta)}|f^{(k)}(z)|^2|g^{(n-k)}(z)|^2d_nA(z)\right)^{q/2}|d\zeta|\\
        & =\int_{\m{T}}\left(\sup_{z\in\Gamma(\zeta)}|g^{(n-k)}(z)|(1-|z|^2)^{n-k}\right)^q\left(\int_{\Gamma(\zeta)}|f^{(k)}(z)|^2dA_k(z)\right)^{q/2}|d\zeta|\\
        & \leq\left(\int_{\m{T}}\sup_{z\in\Gamma(\zeta)}\left(|g^{(n-k)}(z)|(1-|z|^2)^{n-k}\right)^{\frac{pq}{p-q}}|d\zeta|\right)^{\frac{p-q}{p}}\times\\
        & \qquad \times \left(\int_{\m{T}}\left(\int_{\Gamma(\zeta)}|f^{(k)}(z)|^2dA_k(z)\right)^{p/2}|d\zeta|\right)^{q/p}\\
        & \cong \|g\|_{BT^{\frac{pq}{p-q}}}^q\|f\|_{H^p}^q.
    \end{align*}
\end{proof}
\begin{proof}[Proof of Theorem \ref{characterisation for 0<q<p<+infty}] We mention that the suffiency of parts (i) and (ii), are easy conseqeunces of \cite[Theorem 1.3]{chalmoukis2020generalized} and Proposition \ref{Extension of Nikos tgnk result} respectively.
 For the necessity, we start by proving part (ii). Specifically, let $a\in\m{C}^n$ be an arbitrary $n$-tuple, $a\neq \mathbf{0}$, and $T_{g,a}\colon H^p\to H^q$ be bounded. Hence, there exists a constant $C>0$, such that 
$$\|T_{g,a}f\|_{H^q}\leq C\|f\|_{H^p},\qquad \forall f\in H^p.$$
 Fix an $r>0$, $M>1$ and consider an $(R,r)$-lattice $Z=\{z_\lambda\}_{\lambda}$. By  Lemma \ref{lemma of covering hyperbolic disc the stoltz angle}, there exists an $M^*>M>1$ such that
$$\bigcup_{\lambda\colon z_\lambda\in\Gamma_M(\zeta)}D(z_\lambda,r)\subset \Gamma_{M^*}(\zeta).$$
Now, taking into account \eqref{tent norm of Hardy derivatives}, we have that
$$\int_{\m{T}}\left(\int_{\Gamma_{M^*}(\zeta)}\left|\sum_{k=0}^{n-1}f^{(k)}(z)g^{(n-k)}(z)\right|^2dA_n(z)\right)^{q/2}|d\zeta|\leq C\|f\|_{H^p}^q,\quad\forall f\in H^p.$$
Set $$G(z,w)=\mathlarger{\mathlarger{\sum_{k=0}^{n-1}}}\frac{a_k(b)_k|w|^{jk}\overline{w}^kg^{(n-k)}(z)}{(1-|w|^j\overline{w}z)^k}\qquad z,w\in\m{D}\,.$$
 
 Let 
$\mathcal{W}=\{w_\lambda\}_{\lambda}$ be a sequence of points such that
\begin{itemize}
    \item[i)] $w_\lambda\in \overline{D(z_\lambda,r)}$;
    \item[ii)] $w_\lambda$ is a point where the function ${\displaystyle |G(z,z)|(1-|z|^2)^{n}}$ takes its maximum value in $\overline{D(z_\lambda,r)}.$
\end{itemize}
$\mathcal{W}$ may not be hyperbolic separated, as the hyperbolic discs $D(z_\lambda,r)$ overlap. However, as Lemma \ref{finite covering Lemma} implies, we can write it  as a union of finitely many  $(R,r)$-lattices. Hence, without loss of generality, we can assume that $\mathcal{W}$ is an $(R,r)$-lattice.
 Let $b>\max\{1,2/p\}$ and $j\in \m{Z}_{\geq0}.$ Consider as test function
$$S_j[\{c_\lambda\}](z)=\sum_{\lambda}c_\lambda\left(\frac{1-|w_\lambda|^2}{1-|w_\lambda|^j\overline{w_\lambda}z}\right)^b,\qquad \{c_\lambda\}\in T_2^p(\mathcal{W})\,.$$

By  Lemma \ref{Operator Lemma}, we have that
$$\|S_j[\{c_\lambda\}]\|_{H^p}\leq C\|\{c_\lambda\}\|_{T_2^p(\mathcal{W})}.$$

Substituting in the place of $f$ the our test function, we arrive at
$$\int_{\m{T}}\left(\int_{\Gamma_{M^*}(\zeta)}\left|\sum_{\lambda}c_\lambda\left(\frac{1-|w_\lambda|^2}{1-|w_\lambda|^j\overline{w_\lambda}z}\right)^{b}G(z,w_\lambda)\right|^2\,dA_n(z)\right)^{q/2}|d\zeta|\leq C\|c_{\lambda}\|^q_{T_2^p(\mathcal{W})},$$

Let $c_\lambda\in T_2^p(\mathcal{W})$ be a sequence with finetely non-zero entries. Replace $c_\lambda$ in the above formula with $c_\lambda r_\lambda(t)$, where $r_\lambda$ are the Rademacher variables and integrate the resulting inequality with respect of $t$ in $(0,1)$. By applying the \hyperlink{kahane}{Kahane's Inequality}, we arrive at
   \begin{align}
\|c_\lambda\|_{T_2^p(\mathcal{W})}^q&\gtrsim\int_{\m{T}}\int_{0}^{1}\left(\int_{\Gamma_{M^*}(\zeta)}\left|\sum_{\lambda}c_\lambda r_{\lambda}(t)\left(\frac{1-|w_\lambda|^2}{1-|w_\lambda|^j\overline{w_\lambda}z}\right)^{b}G(z,w_\lambda)\right|^2\,dA_n(z)\right)^{q/2}dt\,|d\zeta|\nonumber\\
&\gtrsim\int_{\m{T}}\left(\int_{\Gamma_{M^*}(\zeta)}\int_{0}^{1}\left|\sum_{\lambda}c_\lambda r_{\lambda}(t)\left(\frac{1-|w_\lambda|^2}{1-|w_\lambda|^j\overline{w_\lambda}z}\right)^{b}G(z,w_\lambda)\right|^2\,dt\,dA_n(z)\right)^{q/2}|d\zeta| 
\nonumber\\
&= \int_{\m{T}}\left(\int_{\Gamma_{M^*}(\zeta)}\sum_{\lambda}
|c_\lambda|^2\left(\frac{1-|w_\lambda|^2}{|1-|w_\lambda|^j\overline{w_\lambda}z|}\right)^{2b}d\mu_\lambda(z)\right)^{q/2}|d\zeta|\,,\label{equiKhinichine}
\end{align}
where
$$\mu_{\lambda}(z)= |G(z,w_\lambda)|^2dA_n(z).$$ The last equality holds, due to the fact that $\{r_{\lambda}\}_{\lambda}$ is an orthonormal set in $L^2(0,1).$

As the implicit constant does not depend on the number of points, by a limiting argument, we conclude that \eqref{equiKhinichine} holds also for an arbitrary $\{c_{\lambda}\}_\lambda\in T_2^p(\mathcal{W})\,.$ Using estimates \eqref{geometric estimates}, we readily verify that
$$\chi_{D(w_\lambda,R)}(z)\lesssim\frac{1-|w_\lambda|^2}{|1-|w_\lambda|^j\overline{w_\lambda}z|}\qquad z\in\m{D},$$
where $\chi$ denotes the characteristic function. Consequently, we estimate the right hand side of \eqref{equiKhinichine},
\begin{align*}
    \|c_\lambda\|_{T_2^p(\mathcal{W})}^q&\gtrsim \int_{\m{T}}\left(\int_{\Gamma_{M^*}(\zeta)}\sum_{\lambda}
|c_\lambda|^2\left(\frac{1-|w_\lambda|^2}{|1-|w_\lambda|^j\overline{w_\lambda}z|}\right)^{2b}d\mu_{\lambda}(z)\right)^{q/2}|d\zeta| \\
& \gtrsim \int_{\m{T}}\left(\int_{\Gamma_{M^*}(\zeta)}\sum_{\lambda}
|c_\lambda|^2\chi_{ D(w_\lambda,R)}(z)d\mu_{\lambda}(z)\right)^{q/2}|d\zeta| \\
& \gtrsim \int_{\m{T}}\left(\sum_{\lambda\colon w_\lambda\in \Gamma_M(\zeta)}
|c_\lambda|^2\mu_{\lambda}(D(w_\lambda,R))\right)^{q/2}|d\zeta|\,.
\end{align*}

The last inequality comes from the fact that $\mathcal{W}$ is a finite union of $2R$-hyperbolically separated subsequences.

Let now $\ell>\max\bigl\{1,\frac{1}{q},\frac{p-q}{pq}\bigr\}$ and $\{h_{\lambda}\}_\lambda\in T_1^{\frac{pq\ell}{q-p+p\ell q}}(\mathcal{W})$. By Lemma \ref{factorization of tent sequences} and the fact that $\mathcal{W}$ is a finite union of $(r,R)$-lattices, we can factorize $h_{\lambda}=\Tilde{c}_\lambda t_{\lambda}$ with $\{\Tilde{c}_\lambda\}_\lambda\in T_{2\ell}^{p\ell }(\mathcal{W})$ and $\{t_\lambda\}_\lambda\in T_{\frac{2\ell}{2\ell-1}}^{\frac{\ell q}{\ell q-1}}(\mathcal{W})$, such that
 $$\|\{\Tilde{c}_\lambda\}\|_{T_{2\ell}^{p\ell }(\mathcal{W})}\cdot \|\{t_\lambda\}\|_{T_{\frac{2\ell}{2\ell-1}}^{\frac{\ell q}{\ell q-1}}(\mathcal{W})}\lesssim \|\{h_{\lambda}\}\|_{T_1^{\frac{pq\ell}{q-p+p\ell q}}(\mathcal{W})}\,.$$ 
 
 Set $c_{\lambda}=\Tilde{c}_{\lambda}^{\ell}$. 
 We estimate the following quantity, using first Fubini's Theorem and then H\"older's inequality twice,
\begin{align*}
    \biggl|\sum_{\lambda} & h_\lambda \mu_\lambda^{1/2\ell}(D(w_\lambda,R))(1-|w_\lambda|)\biggr| \\
    &\lesssim   \int_{\m{T}}
    \sum_{\lambda\colon z_\lambda\in\Gamma_M(\zeta)}| t_\lambda \Tilde{c}_\lambda\mu_{\lambda}^{1/2\ell}(D(w_\lambda,R))||d\zeta|  \\
    & \leq \int_{\m{T}}\left(\sum_{\lambda\colon w_\lambda\in\Gamma_M(\zeta)} 
    |c_\lambda|^{2}\mu_{\lambda}\left(D(w_\lambda,R)\right)\right)^{\frac{1}{2\ell}} \left(\sum_{\lambda\colon w_\lambda\in\Gamma_M(\zeta)}|t_\lambda|^{\frac{2\ell}{2\ell-1}}\right)^{1-\frac{1}{2\ell}}|d\zeta|\\
    & \leq \left(\int_{\m{T}}\left(\sum_{\lambda\colon w_\lambda\in\Gamma_M(\zeta)} |c_\lambda|^{2}\mu_{\lambda}\left(D(w_\lambda,R)\right)\right)^{q/2}|d\zeta|\right)^{\frac{1}{q\ell}}\cdot \|\{t_\lambda\}\|_{T_{\frac{2\ell}{2\ell-1}}^{\frac{\ell q}{\ell q-1}}(\mathcal{W})}\\
    & \leq C\|\{\Tilde{c}_\lambda\}\|_{T_2^p(\mathcal{W})}\cdot \|\{t_\lambda\}\|_{T_{\frac{2\ell}{2\ell-1}}^{\frac{\ell q}{\ell q-1}}(\mathcal{W})}\lesssim \|\{h_{\lambda}\}\|_{T_1^{\frac{pq\ell}{q-p+pq\ell}}(\mathcal{W})}
\end{align*}

 By the duality of tent sequences spaces Lemma \ref{duality of tent sequence spaces} and Hahn-Banach we conclude that 
$$\mu_{\lambda}^{1/2\ell}(D(w_\lambda,R))\in T_\infty^{\frac{\ell pq}{p-q}}\iff \mu_{\lambda}(D(w_\lambda,R))\in T_\infty^{\frac{pq}{2(p-q)}},$$
which means that
\begin{equation*}
\int_{\m{T}}\left(\sup_{w_\lambda\in\Gamma_M(\zeta)}\int_{D(w_\lambda,R)}|G(z,w_\lambda)|^2dA_n(z)\right)^{\frac{pq}{2(p-q)}}|d\zeta|<\infty.
\end{equation*}

Using the fact that $|G(z,w_\lambda)|^2$ is subharmonic in the first variable, and estimates \eqref{geometric estimates}, we moreover have that
\begin{align}
\int_{\m{T}}\left(\sup_{w_\lambda \in\Gamma_M(\zeta)}|G(w_\lambda,w_\lambda)|(1-|w_\lambda|^2)^n\right)^{\frac{pq}{p-q}}|d\zeta|<\infty\,.\label{equationfinal}
\end{align}
As before, due to Lemma \ref{lemma of covering hyperbolic disc the stoltz angle}, we can choose 
a $M_{-}$, with $M>M_{-}> 1$  such that
$$\bigcup_{\lambda\colon D(w_{\lambda},R)\cap \Gamma_{M_{-}}(\zeta)\neq\emptyset} D(w_{\lambda}, R) \subset \Gamma_{M}(\zeta)$$

Finally,
\begin{align}
    \int_{\m{T}}&\left(\sup_{z\in\Gamma_{M_-}(\zeta)}\left|\sum_{k=0}^{n-1}(b)_ka_k|z|^{jk}\overline{z}^kg^{(n-k)}(z)\frac{(1-|z|^2)^n}{(1-|z|^{j+2})^{k}}\right|\right)^{\frac{pq}{p-q}}|d\zeta|= \nonumber\\
    &=\int_{\m{T}}\left(\sup_{z\in\Gamma_{M_{-}}(\zeta)}|G(z,z)|(1-|z|^2)^n\right)^{\frac{pq}{p-q}}|d\zeta|\nonumber\\
    & \leq \int_{\m{T}}\left(\sup_{w_\lambda\in\Gamma_{M}(\zeta)}|G(w_{\lambda},w_{\lambda})|(1-|w_{\lambda}|^2)^n\right)^{\frac{pq}{p-q}}|d\zeta|\nonumber<\infty
\end{align}
due to \eqref{equationfinal}. The next step is to use Lemma \ref{algebraic lhmma} for $f_k(z)=a_k(b)_k\overline{z}^kg^{(n-k)}(z)$. In the notation of Lemma \ref{algebraic lhmma} we have proved that 
$$D_j\in T^{\frac{pq}{p-q}}_{\infty}\qquad j=0,\dots,n-1.$$

 Therefore, for any $k$ such that  $ a_k\neq 0$, the function $g^{(n-k)}(z)(1-|z|)^{n-k}$ in $ \frac{1}{2}<|z|<1$ can be written as a linear combination of products of bounded functions and the functions $D_j$. Hence,
$$\int_{\m{T}}\left(\sup_{z\in\Gamma(\zeta)}|g^{(n-k)}(z)|(1-|z|)^{n-k}\right)^p|d\zeta|<\infty.$$

Proposition \ref{Blochtentderivatives} now implies that
$$g\in BT^{\frac{pq}{p-q}}.$$
 So far, we have proved that if $T_{g,a}$ is bounded, then $g\in BT^{\frac{pq}{p-q}}$, which proves part (ii) of the theorem.

 To prove part (i), we recall that, since $g\in BT^{\frac{pq}{p-q}}$, Proposition \ref{Extension of Nikos tgnk result} implies that $T_g^{n,k}\colon H^p\to H^q$ are bounded for $1\leq k \leq n-1.$ Therefore, when $a_0\neq 0$, the boundedness of $T_{g,a}$ implies the boundedness of $T_{g}^{n,0}$, since
$$a_0T_g^{n,0}=T_{g,a}-\sum_{k=1}^{n-1}a_kT_g^{n,k}.$$
By an application of \eqref{tent norm of Hardy derivatives}, the boundedness of $T_{g}^{n,0}$ is equivalent to
    $$\int_{\m{T}}\left(\int_{\Gamma(\zeta)}|f(z)|^2|g^{(n)}(z)|^2dA_n(z)\right)^{q/2}|d\zeta|\leq C\|f\|^q_{p}\qquad\forall f\in H^p.$$
   However, \cite[Theorem 1.2]{lv2021pautent} implies that the above inequality is equivalent to the fact that
    $$\int_{\m{T}}\left(\int_{\Gamma(\zeta)}|g^{(n)}(z)|^2\,dA_n(z)\right)^{\frac{pq}{2(p-q)}}|d\zeta|<\infty.$$
    Finally, using once more \eqref{tent norm of Hardy derivatives}, we acquire that $g\in H^{\frac{pq}{p-q}}.$\\
\end{proof}

\section{Compactness of \texorpdfstring{$T_{g,a}$}\empty}

The compactness of $T_{g,a}$ relies  on the compactness of the Volterra operator $V$ on Hardy spaces and an approximation argument. We recall from \cite[Lemma 1]{aleman2001integral} that $V\colon H^p\to H^q$ is compact, for $0<p,q<\infty$ with $\frac{1}{p}-\frac{1}{q}<1$. Consequently, we prove the following auxiliary proposition.
\begin{prop}\label{Prop for polynomil Tp compactness}
     Let $0<p<q<\infty$, $n\in\m{N}$, $0\leq k\leq n-1$ fixed and $g\in \Hol(\m{D})$. Set $\alpha=\frac{1}{p}-\frac{1}{q}$. If $\kappa<\alpha\leq \kappa+1< n-k$ for some $\kappa\in\m{Z}_{\geq0}$ and $ g^{(\kappa)}\in\lambda_{\alpha-k}$, then $T_g^{n,k}\colon H^q\to H^q$ is compact.
    
 \end{prop}
 \begin{proof}
        We start by proving that if $r\in\m{N}$ and $\alpha=\frac{1}{p}-\frac{1}{q}<r$ and $p<\frac{1}{r-1}$, then $V^{r}\colon H^p\to H^q$ is compact. For, if $r=1$, then this reduces to \cite[Lemma 1]{aleman2001integral}. If $r>1$, then we write $V^r=VV^{r-1}$ and we verify, by using \cite[Theorem 1]{aleman2001integral}, that  $V^{r-1}\colon H^p\to H^s$ is bounded, where $s=\frac{p}{1-(r-1)p}$\,. Then $\frac{1}{s}-\frac{1}{q}<1$ and consequently $V\colon H^s\rightarrow H^q$ is compact. Hence, the operator $V^{r}$ is also compact. 
        
        Notice then, that for $r=n-k$, our assumption implies that $p<\frac{1}{n-k-1}$ and therefore, $V^{n-k}\colon H^p\to H^q$ is compact.
        Now, let $P$ be a polynomial. 
       Consecutive integration by parts, one can write $T_{P}^{n,k}$ in the following form
       \begin{align*}
           T_{P}^{n,k}&=c_1V^{n-k} M_{P^{(n-k)}}+c_2 V^{n-k+1} M_{P^{(n-k+1)}}+\cdots+c_kV^{n}M_{P^{(n)}}\\
           &=V^{n-k} (c_1 M_{P^{(n-k)}}+c_2  V M_{P^{(n-k+1)}}+\cdots+c_k V^kM_{P^{(n)}})
       \end{align*}
       where $M_{P^{(j)}}f=f\cdot P^{(j)}$ is the pointwise multiplication operator. Moreover, the multiplication operators induced by polynomial symbols are bounded on $H^p$, and so are the operators $V, \dots, V^k$. Thus, we have that $T_P^{n,k}$ is compact.
        Let now $g\in \Hol(\m{D})$ satisfying $g^{(\kappa)}\in\lambda_{\alpha-k}$. 
       The polynomials are dense in $\lambda_{\alpha-\kappa}$, (see for example  \cite[Proposition 2]{zhu1993bloch}), hence there exists a sequence of polynomials $\{P_m\}_{m\in\m{N}}$ such that
       $\|g^{(\kappa)}-P_m\|_{\Lambda_{\alpha-\kappa}}\to 0$ as $m$ goes to infinity. Consider the polynomials satisfying 
       $$\begin{cases}
     G_{m}^{(\kappa)}(z)=P_m(z), \,\,\, z\in\m{D}\\
     G_{m}(0)=\dots=G_{m}^{(\kappa-1)}(0)=0.
 \end{cases}$$
       Then for the operator norms, we have that
       \begin{align*}
           \|T_{g}^{n,k}-T_{G_m}^{n,k}\|=\|T_{g-G_m}^{n,k}\|\leq C\|g^{(\kappa)}-P_m\|_{\Lambda_{\alpha-\kappa}}\rightarrow 0\,,\quad m\to\infty\,.
       \end{align*}
       As $T_{g}^{n,k}$ is approximated in operator norm by a sequence of compact operators, it is compact.\\
 \end{proof}
\begin{proof}[Proof of Theorem \ref{compactness characterisation on Hardy spaces}](i)

Assume that $\alpha<n-\ell$ and $g$ satisfies (\ref{equation for p<q in Tga compact Hardy}). Proposition \ref{Proposition Lipschitz class derivatives} implies this is equivalent to
$$ \lim_{|z|\to 1}\frac{(1-|z|)^{n-k}|g^{(n-k)}(z)|}{(1-|z|)^{\alpha}}=0.\qquad \forall \,0\leq k\leq \ell.$$

Hence, Proposition \ref{Prop for polynomil Tp compactness} implies that, $T_g^{n,k}$ are simultaneously compact for $0\leq k\leq \ell$ and consequently $T_{g,a}$ is compact as a sum of compact operators.
 For the other implication, let $\{\lambda_n\}_n\subset\m{D}$ such that $|\lambda_n|\to 1$ as $n\to\infty$. Use the test functions
 $$f_{\lambda_n,\gamma}(z)= \frac{(1-|\lambda_n|^2)^{\gamma-1/p}}{(1-\overline{\lambda_n}z)^\gamma}\qquad \gamma>1/p,\,\,z\in\m{D}.$$

These functions converge uniformly to zero on compact subsets of the unit disc as $n\to \infty$ and there exists a constant $C=C(\gamma)>0$ such that $\|f_{\lambda,\gamma}\|_{H^p}\leq C$.

 So, the growth estimate \eqref{equation for pointevaluations} shows that
         $$|T_{g,a}(f_{\lambda_n,\gamma})^{(n)}(\lambda_n)|\lesssim\frac{\|T_{g,a}(f_{\lambda_n})\|_{H^q}}{(1-|\lambda_n|^2)^{1/q+n}}.$$
      The compactness of $T_{g,a}$ implies that
         $\|T_{g,a}(f_{\lambda_n,\gamma})\|_{H^q}\rightarrow0$ as $n$ goes to infinity. Computing the $n$-th derivative of $T_{g,a}$, we verify that
      $$\lim_{n\to\infty}\left|\sum_{k=0}^{n-1}a_k\overline{\lambda_n}^k(\gamma)_k\frac{(1-|\lambda_n|^2)^{n-k}g^{(n-k)}(\lambda_n)}{(1-|\lambda_n|^2)^{\alpha}}\right|=0\,.$$
      Appealing to Lemma \ref{Nikos Lemma for test functions}, we conclude that
         $$\lim_{n\to\infty}\frac{(1-|\lambda_n|^2)^{n-\ell}|g^{(n-\ell)}(\lambda_n)|}{(1-|\lambda_n|^2)^{\alpha}}=0$$
         which is equivalent to $g^{(\kappa)}\in\lambda_{\alpha-\kappa}\,.$
         
          For (b), we observe that similar steps as the in the previous part, allow us to conclude that $\lim_{n\to\infty}|g^{(n-\ell)}(\lambda_n)|=0,$
         which implies that $g^{(n-\ell)}\equiv 0$, hence $T_{g,a}$ is the zero operator.
          For (ii), let $P$ be a polynomial and $0\leq k\leq n-1$. As a consequence of \cite[Proposition 3.5]{chalmoukis2020generalized}, $T_P^{n,k}$ acts compactly on $H^p$, therefore  $T_P^{n,k}\colon H^p\to H^q$ is compact too. 
         The arguments now follow the same pattern as in (i). Specifically let $g\in H^{\frac{pq}{p-q}}$. As polynomials are dense in Hardy spaces, there exists a sequence of polynomials approximating $g$ in $\|\cdot \|_{H^{\frac{pq}{p-q}}}$ and consequently $T_{g}^{n,k}$ is the limit in operator norm of compact operators, therefore compact. So $T_{g,a}$ is compact, as linear combination of compact operators. \\
\end{proof}

\bibliographystyle{plain}
\bibliography{bibliography}

\end{document}